\documentclass{amsart}
\usepackage[utf8x]{inputenc}
\usepackage{amssymb,amsmath}
\usepackage{setspace}
\usepackage[dvips]{graphicx}

\newtheorem{theo}{Theorem}[section]

\newtheorem{coro}{Corollary}[section]

\def\Ric{{\rm Ric}}
\def\Ricf{{\rm Ric}^f}

\def\ricf{{\rm ric}^f}

\title{On the Geometry of Warped Foliations}
\author{Szymon M. Walczak}
\date{version \today}

\begin{document}

\begin{abstract}
We discuss the geometry of warped foliations. After examining the Levi-Civita connection, and the curvature tensor, we describe the formulae for sectional, Ricci and scalar curvatures. In the final part of this note, we present some examples. 
\end{abstract}

\maketitle

\section{Introduction}

The notion of the warped foliation was introduced in \cite{W4} by the author of this note. It is a generalization of the M. Berger's modification of a Riemannian structure of $S^3$ along the fibers of the Hopf fibration called Berger shperes \cite{P}.

Warped foliations were widely studied in the point of view of the Gromov-Hausdorff convergence \cite{G}. The sufficient and necessary condition of converging in Gromov–Hausdorff sense of a Riemannian submersion  and Riemannian foliation with all leaves compact to the space of leaves with a metric defined by Hausdorff distance of leaves were already developed in \cite{W1} and \cite{W2}. Moreover, in \cite{W3}, the author of this paper has presented the connection betweeen the Hausdorff leaf space for a given foliation and the warped foliation.

\section{Preliminaries}

At the beginning, let us recall the definitions needed in this note.

Let $(M,g)$ be a Riemannian manifold, while $\nabla$ the Levi-Civita connection on $M$. Let $P$ be a smooth distribution on $M$. Following \cite{R}, one can define a smooth tensor fields $T$ of type $(1,2)$ by the formula
\begin{eqnarray}
g(T(U,V), X) &=& g(\frac{1}{2} \nabla_U V + \frac{1}{2}\nabla_V U,X),
\nonumber\\
g(T(U,X), V) &=& -g(T(U,V), X),
\nonumber\\
g(T(X,\cdot), Y) &=& g(T(X,\cdot), U) = 0,
\nonumber
\end{eqnarray}
and  $A$ (also of type $(1,2)$) by
\begin{eqnarray}
g(A(U,V), X) &=& g(\frac{1}{2} \nabla_U V - \frac{1}{2}\nabla_V U,X),
\nonumber\\
g(A(U,X), V) &=& -g(A(U,V), X),
\nonumber\\
g(A(X,\cdot), Y) &=& g(A(X,\cdot), U) = 0,
\nonumber
\end{eqnarray}
where $U,V\in P$, and $X,Y\in P^{\bot}$. $T$ is called \textit{the second fundamental form of $P$}, while $A$ \textit{the integrability tensor}. It follows from Frobenius Theorem \cite{KN} that for $P$ integrable the integrability tensor $A$ vanishes.

Now, consider a foliation\footnote{ Theory of foliations can be found in \cite{CC}, the amazing book written by A. Candel and L. Conlon} $\mathcal{F}$ on the manifold $(M,g)$. There are two natural distributions on $(M,g)$ defined by $\mathcal{F}$. One of them, consisting of all vectors tangent to the leaves of $\mathcal{F}$ is called {\em tangent}, and will be also denoted by $\mathcal{F}$. The second one, called {\em  orthogonal} consists of all vectors which are $g$-orthogonal to the leaves of $\mathcal{F}$. It will be denoted by $\mathcal{F}^{\bot}$. 

Denote by $T$ the second fundamental form of $\mathcal{F}$. Since $\mathcal{F}$ is integrable, its integrability tensor vanishes everywhere. Let $S$ denotes the second fundamental form, while $A$ the integrability tensor of $\mathcal{F}^{\bot}$. 

Let us recall that a foliation $\mathcal{F}$ satisfying 
\begin{equation}\label{eq:Riemannian foliation}
 \mathcal{L}_U g(X,Y) = 0,
\end{equation}
where $\mathcal{L}$ denotes the Lie differentiation, $U\in\mathcal{F}$, $X,Y\in\mathcal{F}^{\perp}$, is called {\em Riemannian foliation}\cite{MM}.

\begin{theo}\label{thm:Riemannain foliation and the second fundamental form}
 $\mathcal{F}$ is Riemannian if and only if the second fundamental form $S$ of the orthogonal distribution $\mathcal{F}^{\bot}$ vanishes.
\end{theo}

A proof can be found in \cite{R}.

Throughout this paper, we will also use the following denotations:
\[
H_f(U) = \nabla_U \nabla f,\quad h_f (X,Y) = \langle \nabla_X \nabla f, Y\rangle.
\]

\section{Warped foliations}

Let $(M,\mathcal{F},g)$ be a foliated Riemannian manifold. Consider a smooth function $f:M\to (0,\infty)$ constant along the leaves of $\mathcal{F}$. Such function is often called a basic function. Modify the Riemannian structure $g$ to $g_f$ in the following way: 

Let $g_f (v,w) = f^2 g(v,w)$ for vectors $v,w$ tangent to the foliation $\mathcal{F}$. Next, if at least one of vectors $v,w$ is perpendicular to $\mathcal{F}$ then set $g_f (v,w) = g(v,w)$. Foliated Riemannian manifold $(M,\mathcal{F},g_f)$ we call \textit{the warped foliation} and denote by $M_f$. The function $f$ is called \textit{the warping function}.

\begin{figure}[h]
 \centering
 \includegraphics{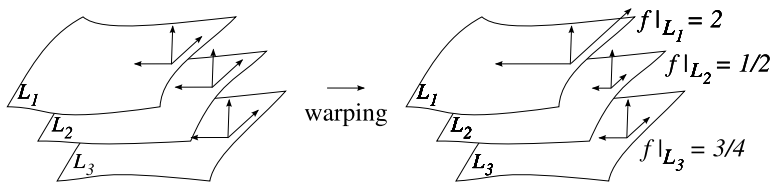}
 \caption{Warping of a foliation.}
 \label{fig:Warping of a foliation}
\end{figure}

Briefly speeking, one can understand warping as a conformal modification of a Riemannian structure of a foliated Riemannian manifold $(M,\mathcal{F},g)$ along the leaves of $\mathcal{F}$ with a dilatation equal, for given leaf, to the value of a warping function on this leaf (see Fig. \ref{fig:Warping of a foliation}). The orthogonal vectors remain unchanged.

In the following section, we will study the geometry of the warped foliations on a compact Riemannian foliated manifolds. We will calculate the Levi-Civita connection, the curvature tensor, and the curvatures of the warped foliations.

\section{Levi-Civita connection and curvature tensor}\label{sec:Levi-Civita connection and curvature tensor}

Let $(M,\mathcal{F},g)$ be a foliated Riemannian manifold of dimension $p$ and codimension $q$. Let $\nabla$ deonte the Levi-Civita connection on $(M,\mathcal{F},g)$, and let $f:M\to (0,\infty)$ be a warping function on $M$. Denote by $\nabla^f$ the Levi-Civita connection of the warped foliation $M_f$.

\begin{theo}\label{thm:Levi-Civita connection and curvature tensor}
The Levi-Civita $\nabla^f$ on $M_f$ is given by
\[
\nabla^f _X Y = (\nabla _X Y)^{\bot} + \frac{1}{f^2} (\nabla _X Y)^{\top} - \frac{1-f^2}{f^2} A(X,Y),
\]
\[
\nabla^f _U V = (\nabla _U V)^{\top} + f^2 (\nabla _U V)^{\bot} - \frac{1}{2}\langle U, V\rangle\cdot \nabla f^2,
\]
\[
\nabla^f _X U = \nabla _X U + \frac{1}{2} \frac{X f^2}{f^2} U - (1-f^2) A(X,U),
\]
\[
\nabla^f _U X = \nabla _U X + \frac{1}{2} \frac{X f^2}{f^2} U - (1-f^2) A(X,U),
\]
where $X$ and $Y$ are ortoghonal, but $U$ and $V$ tangent to $\mathcal{F}$.
\end{theo}
\begin{proof}
By the definition of $g_f$, and by the Koszul formula
\begin{align}
g_f(\nabla_E F,G) &= \frac{1}{2}( Eg_f(F,G) + Fg_f(G,E) - Gg_f(E,F)
\nonumber\\
& + g_f([E,F],G) - g_f([F,G],E) + g_f([G,E],F)).
\nonumber
\end{align}
Direct calculations give the statement.
\end{proof}

Following Theorem \ref{thm:Levi-Civita connection and curvature tensor}, one can calculate the curvature tensor $R^f$ for the warped foliation $M_f$. Let $W,X,Y,Z\in\mathcal{F}^{\bot}$ and $P,Q,U,V\in\mathcal{F}$ be vector fields on $M$. We have
\begin{align}
\langle R^f(X,Y)Z,W\rangle_f &= \langle R(X,Y)Z,W\rangle - 2(1-f^2)\langle A(X,Y), A(Z,W)\rangle
\nonumber\\
&+ \langle A(X,W), A(Y,Z)\rangle - \langle A(X,Z), A(Y,W)\rangle
\nonumber\\
&- \frac{1-f^2}{f^2} (\langle S(X,W), S(Y,Z)\rangle - \langle S(X,Z),S(Y,W)\rangle),
\nonumber
\end{align}
\begin{align}
\langle R^f(X,Y)Z,U \rangle_f &= \langle R(X,Y) Z, U\rangle + f\cdot Xf\cdot  \langle A(Y,Z), U \rangle 
\nonumber\\
&- f\cdot Yf\cdot  \langle A(X,Z), U \rangle -2f\cdot Zf\cdot  \langle A(X,Y), U \rangle
\nonumber\\
&- \frac{Xf}{f}\langle S(Y,Z),U \rangle - \frac{Yf}{f}\langle S(X,Z),U \rangle
\nonumber\\
&- (1-f^2)(\langle(\nabla_X A)(Y,Z), U \rangle + \langle(\nabla_Y A)(X,Z), U \rangle
\nonumber\\
&- \langle A(X,Y),T(U,Z) \rangle),
\nonumber
\end{align}
and
\begin{align}
\langle R^f(X,U)Y,V \rangle_f &= f^2 \langle R(X,U)Y,V\rangle + f \langle \nabla_X \nabla f, Y\rangle \langle U,V\rangle
\nonumber\\
&- f\cdot Xf\cdot \langle T(U,V), Y\rangle - f\cdot Yf\cdot  \langle T(U,V), X \rangle 
\nonumber\\
&+ f^2 (1-f^2) \langle A(X,V), A(Y,U) \rangle
\nonumber\\
&+ (1-f^2) [\langle S(X,V), A(Y,U) \rangle
\nonumber\\
&- \langle S(Y,V), (\nabla_X U)^{\bot} \rangle - \langle (\nabla _U S) (X,Y),V \rangle]
\nonumber
\end{align}
Moreover, since
\[
\langle R^f(U,V)X,Y\rangle + \langle R^f(V,X)U,Y\rangle + \langle R^f(X,U)V,Y\rangle = 0,
\]
then
\begin{align}
\langle R^f(U,V)X,Y \rangle_f &= f^2 \langle R(U,V)X,Y\rangle
\nonumber\\
&+ f^2(1-f^2) [\langle A(X,V),A(Y,U)\rangle - \langle A(Y,V),A(X,U)\rangle]
\nonumber\\
&+ 2(1-f^2)[\langle S(X,V),A(Y,U)\rangle - \langle S(Y,V),A(X,U)\rangle].
\nonumber
\end{align}
In addition,
\begin{align}
\langle R^f(U,V)P,X \rangle_f &= f^2 \langle R(U,V)P,X\rangle
\nonumber\\
&+ f^2(1-f^2)\langle A(X,U),T(V,P)\rangle 
\nonumber\\
&- f^2(1-f^2) \langle A(X,V),T(U,P)\rangle 
\nonumber\\
&- f(1-f^2)\langle V,P\rangle \langle A(\nabla f, X), U\rangle
\nonumber\\
&+ f(1-f^2)\langle U,P\rangle \langle A(\nabla f, X), V\rangle
\nonumber\\
&+ f\langle U,P\rangle \langle H_f(V),X\rangle - f\langle V,P\rangle \langle H_f(U),X\rangle,
\nonumber
\end{align}
and
\begin{align}
\langle R^f(U,V)P,Q\rangle_f &= f^2\langle R^{\top} (U,V)P,Q\rangle
\nonumber\\
&- f^3\langle U,P\rangle \langle \nabla f, T(V,Q\rangle - f^3\langle V,Q\rangle \langle \nabla f, T(U,P)\rangle
\nonumber\\
&+ f^4\langle T(U,P),T(V,Q\rangle + f^2 \langle U,P\rangle \langle V,Q\rangle \|\nabla f\|^2
\nonumber\\
&+ f^3\langle V,P\rangle \langle \nabla f, T(U,Q\rangle + f^3\langle U,Q\rangle \langle \nabla f, T(V,P)\rangle
\nonumber\\
&- f^4\langle T(V,P),T(U,Q\rangle - f^2 \langle V,P\rangle \langle U,Q\rangle \|\nabla f\|^2.
\nonumber
\end{align}
where $R$ denotes the curvature tensor of $(M,g)$.

\section{Curvatures of warped foliations}\label{sec:Curvatures of warped foliations}

We are now able to calculate the curvatures for $M_f$. Let $U$, $V$ be tangent, while $X$, $Y$ orthogonal vectors tangent to $M_f$ at point $x$. The natural consequence of the above section is the following theorem.

\begin{theo}\label{thm:sectional curvature}
The sectional curvature $\kappa^f$ of a warped foliations $M_f$ satisfies
\begin{align}
\kappa ^f (X,Y) &= \kappa(X,Y) + 3(1-f^2)\|A(X,Y)\|^2
\nonumber\\
&+ \frac{1-f^2}{f^2}\|S(X,Y)\|^2 - \frac{1-f^2}{f^2}\langle 
S(X,X),S(Y,Y) \rangle,
\nonumber\\
\kappa ^f (X,U) &= \kappa (X,U) - \frac{1}{f} h_f(X,X) + 2fXf \langle 
T(U,U), X\rangle
\nonumber\\
&- (1-f^2) [\langle (\nabla_U S)(X,X), U\rangle - \|S(X,U)\|^2]
\nonumber\\
&- f^2(1-f^2) \|A(X,U)\|^2,
\nonumber\\
\kappa ^f (U,V) &= \frac{\hat \kappa(U,V)}{f^2} - \frac{\|\nabla 
f\|^2}{f^2} - f^4 \langle T(U,U), T(V,V)\rangle 
\nonumber\\
&+ f^4 \|T(U,V)\|^2 + f\langle \nabla f, T(V,V)\rangle + f\langle 
\nabla f, T(U,U)\rangle,
\nonumber
\end{align}
where $\hat \kappa$ denotes the sectional curvature of a leaf.
\end{theo}
\begin{proof}
Follows directly from the formulae for curvature tensor from the previous section.
\end{proof}

We will now calculate the Ricci tensor, and the Ricci curvature for the warped foliation $M_f$. Let $U_1,\dots, U_p,X_1,\dots, X_q$ be an orthogonal basis on $M_f$ in a point $x$. Let us recall, that for any $E,F\in T_x M_f$ we have
\[
\Ricf (E,F) = \sum _{i=1} ^{p} \langle R^f(U_i,E)F, U_i\rangle_f + \sum _{j=1} ^{q} \langle R^f(X_j, E)F, X_j\rangle_f,
\]
where $\langle\cdot,\cdot\rangle_f = g_f(\cdot,\cdot)$. Set $\bar U_i= fU_i$, $i=1,\dots,p$. $\bar U_1,\dots,\bar U_p,X_1,\dots,X_q$ form orthogonal basis in $T_x M$. By the results of Section \ref{sec:Levi-Civita connection and curvature tensor}, for any $U,V$ tangent to $\mathcal{F}$
\begin{align}
\langle R^f(U_i,U)V,U_i\rangle _f &= \langle R^{\top}(\bar U_i,U)V,\bar U_i\rangle 
\nonumber\\
&- f^2\langle T(\bar U_i,\bar U_i),T(U,V)\rangle + \langle T(\bar U_i,U),T(\bar U_i,V)\rangle
\nonumber\\
&- \langle U,V\rangle \|\nabla f\|^2 + \langle\bar U_i,U\rangle \langle\bar U_i, V\rangle \|\nabla f\|^2
\nonumber\\
&+ f\langle T(U,V),\nabla f\rangle + f\langle U,V\rangle \langle T(\bar U_i,\bar U_i),\nabla f\rangle
\nonumber\\
&+ f\langle \bar U_i,U\rangle \langle T(\bar U_i,V),\nabla f\rangle + f\langle \bar U_i,V\rangle \langle T(\bar U_i,U),\nabla f\rangle
\nonumber
\end{align}
and
\begin{align}
\langle R^f(X_j,U)V,X_j\rangle _f &= f^2 \langle R(X_j,U)V,X_j\rangle
\nonumber\\
&+ (1-f^2)[\langle S(X_j,U),S(X_j,V)\rangle+\langle \nabla_U S(X_j,X_j),V\rangle]
\nonumber\\
&+ f^2(1-f^2) \langle A(X_j,U),A(X_j,V)\rangle 
\nonumber\\
&- 2fX_j f\langle T(U,V),X_j\rangle
\nonumber\\
&+ f\langle \nabla_{X_j}\nabla f,X_j \rangle\langle U,V\rangle.
\nonumber
\end{align}
Recall that
\[
\Ricf (U,V) = \sum_{i=1}^{p} \langle R^f(U_i,U)V,U_i\rangle _f + \sum_{i=1}^{q} \langle R^f(X_j,U)V,X_j\rangle _f.
\]
Finally
\begin{align}\label{eqn:Ricci curvature for two tangent vectors}
\Ricf (U,V) &= \Ric ^{\mathcal{F}} (U,V) + f^2\Ric ^{\bot} (U,V)
\\
&- f^2 \langle H^{\mathcal{F}}, T(U,V) \rangle + f\langle U,V\rangle\langle H^{\mathcal{F}},\nabla f \rangle
\nonumber\\
&- (p-1) \langle U,V\rangle \|\nabla f\|^2 + pf\langle T(U,V),\nabla f\rangle
\nonumber\\
&+ f^2 \langle T^{\top} U,T^{\top} V\rangle - f^2(1-f^2) \langle A^{\bot}U,A^{\bot} V\rangle
\nonumber\\
&+ (1-f^2) \langle S^{\bot} U, S^{\bot} V\rangle + (1-f^2) {\rm tr} ^{\bot}\langle (\nabla_U S)(\cdot,\cdot), V\rangle
\nonumber\\
&- f\langle U,V\rangle {\rm tr}^{\bot} h_f
\nonumber
\end{align}
where
\begin{align}
\langle T^{\top} U,T^{\top} V\rangle &= \sum_{i=1}^{p} \langle T_{U_i} U,T_{U_i} V\rangle,
\nonumber\\
\langle A^{\bot} U,A^{\bot} V\rangle &= \sum_{i=1}^{q} \langle A_{X_i} U,A_{X_i} V\rangle,
\nonumber\\
\langle S^{\bot} U,S^{\bot} V\rangle &= \sum_{i=1}^{q} \langle S_{X_i} U,S_{X_i} V\rangle,
\nonumber\\
H^{\mathcal{F}} &= \sum_{i=1}^{p} T(U_i,U_i),
\nonumber\\
tr^{\bot} F(\cdot,\cdot) &= \sum_{i=1}^{q} F(X_i,X_i).
\nonumber
\end{align}

Now, let $X,Y$ be orthogonal to $\mathcal{F}$. Similarly, we get
\begin{align}\label{eqn:Ricci curvature for two orthogonal vectors}
\Ricf (X,Y) &= \Ric^{\bot}(X,Y) + \Ric^{\top}(X,Y)
\\
&- (1-f^2)\langle A_X^{\top},A_Y^{\top}\rangle + 3(1-f^2)\langle A^{\bot} X, A^{\bot} Y\rangle
\nonumber\\
&+ \frac{Xf}{f}\langle H^{\mathcal{F}},Y\rangle + \frac{Yf}{f}\langle H^{\mathcal{F}},X\rangle
\nonumber\\
&- p\frac{h_f(X,Y)}{f} + \frac{1-f^2}{f^2}\langle H^{\bot},S(X,Y)\rangle
\nonumber\\
&- \frac{1-f^2}{f^2}[\langle S_X^{\top},A_Y^{\top}\rangle - {\rm tr}^{\top} \langle(\nabla_{\cdot} S)(X,Y),\cdot \rangle
\nonumber\\
&- \langle S_Y^{\top},(\nabla_X^{\bot})^{\top}\rangle - \langle S_X^{\bot},S_Y^{\bot}\rangle],
\nonumber
\end{align}
with
\begin{align}
H^{\bot} &= \sum_{i=1}^{q} T(X_i,X_i),
\nonumber\\
\langle S_X^{\top},A_Y^{\top}\rangle &= \sum_{i=1}^{p} \langle S(X,U_i), A(Y,U_i)\rangle,
\nonumber\\
{\rm tr}^{\top} \langle(\nabla_{\cdot} S)(X,Y),\cdot \rangle &= \sum_{i-1}^{p}\langle(\nabla_{U_i} S)(X,Y),U_i \rangle,
\nonumber\\
\langle S_Y^{\top},(\nabla_X^{\bot})^{\top}\rangle &= \sum_{i=1}^{q} \langle S(Y,U_i),\nabla_X^{\bot} U_i\rangle,
\nonumber\\
\langle S_X^{\bot},S_Y^{\bot}\rangle &= \sum_{i=1}{q}\langle S(X,X_i),S(Y,X_i)\rangle.
\nonumber
\end{align}

Finally, if $X$ is orthogonal and $U$ tangent to $\mathcal{F}$, we have
\begin{align}\label{eqn:Ricci curvature for orthogonal and tangent vectors}
\Ricf (X,U) &= \Ric^{\bot}(X,U) + \Ric^{\top}(X,U)
\\
&+ 3f\langle A(\nabla f,X),U\rangle + (p-1)\frac{(1-f^2)}{f}\langle A(\nabla f,U),X\rangle
\nonumber\\
&- \frac{p-1}{f} h_f(U,X) + (1-f^2) \langle A_X^{\top},T_U^{\top}\rangle
\nonumber\\
&+ (1-f^2)\langle A(U,X),H^{\mathcal{F}}\rangle + \frac{1}{f}\langle S(X,\nabla f), U\rangle 
\nonumber\\
&+ (1-f^2) {\rm tr}^{\bot} \langle(\nabla _{\cdot} A)(X,\cdot) - (\nabla_X A)(\cdot,\cdot),U\rangle
\nonumber\\
&- \frac{Xf}{f}\langle H^{\bot}, U\rangle,
\nonumber
\end{align}
where $(\nabla _{\cdot} A)(X,\cdot) = \sum_{i=1}^{q}(\nabla_{X_i} A)(X,X_i).$

Let $E\in T_x M_f$ be an unit vector. We have $E=aU + bX$, where $|U|_f=|X|_f=1$, $U\in\mathcal{F}^{\top}$, $X\in\mathcal{F}^{\bot}$ and $a^2+b^2=1$. Note that Ricci curvature $\ricf$ in a point $x$ for the warped foliation $M_f$ is given by the formula
\begin{eqnarray}
\ricf (E) &=& \Ricf (E,E)
\nonumber\\
&=& a^2\Ricf (U,U) + 2ab \Ricf (X,U) + b^2 \Ricf (X,X).
\nonumber
\end{eqnarray}

\begin{theo}\label{thm:Ricci curvature}
The Ricci curvature $\ricf$ in a point $x$ for the warped foliation $M_f$ satisfies
\begin{align}
\ricf(E) &= a^2(\Ric^{\mathcal{F}} (U,U) + f^2\Ric ^{\bot} (U,U)) +2ab(\Ric (X,U) + b^2\Ric (X,X)
\nonumber\\
&+ a^2(-f^2 \langle H^{\mathcal{F}}, T(U,U) \rangle + \frac{1}{f}\langle H^{\mathcal{F}},\nabla f \rangle - \frac{(p-1)\|\nabla f\|^2}{f^2} 
\nonumber\\
&+ pf\langle T(U,U),\nabla f\rangle + f^2 \|T^{\top} U\|^2 - f^2(1-f^2)\|A^{\bot}U\|^2
\nonumber\\
&+ (1-f^2) \|S^{\bot} U\|^2 + (1-f^2) {\rm tr} ^{\bot}\langle (\nabla_U S)(\cdot,\cdot), U\rangle - \frac{{\rm tr}^{\bot} h_f}{f})
\nonumber\\
&+ 2ab( 3f\langle A(\nabla f,X),U\rangle + (p-1)\frac{(1-f^2)}{f}\langle A(\nabla f,U),X\rangle
\nonumber\\
&- \frac{p-1}{f} h_f(U,X) + (1-f^2) \langle A_X^{\top},T_U^{\top}\rangle + (1-f^2)\langle A(U,X),H^{\mathcal{F}}\rangle 
\nonumber\\
&+ \frac{1}{f}\langle S(X,\nabla f), U\rangle + (1-f^2) {\rm tr}^{\bot} \langle(\nabla _{\cdot} A)(X,\cdot) - (\nabla_X A)(\cdot,\cdot),U\rangle
\nonumber\\
&- \frac{Xf}{f}\langle H^{\bot}, U\rangle)
\nonumber\\
&+ b^2(-(1-f^2)\|A_X^{\top}\|^2 + 3(1-f^2)\|A^{\bot} X\|^2 + \frac{2Xf}{f}\langle H^{\mathcal{F}},X\rangle 
\nonumber\\
&- p\frac{h_f(X,X)}{f} + \frac{1-f^2}{f^2}\langle H^{\bot},S(X,X)\rangle - \frac{1-f^2}{f^2}[\langle S_X^{\top},A_X^{\top}\rangle 
\nonumber\\
&- {\rm tr}^{\top} \langle(\nabla_{\cdot} S)(X,X),\cdot \rangle - \langle S_X^{\top},(\nabla_X^{\bot})^{\top}\rangle - \langle S_X^{\bot},S_X^{\bot}\rangle]).
\nonumber
\end{align}
\end{theo}
\begin{proof}
 The proof follows directly from the formulae (\ref{eqn:Ricci curvature for two tangent vectors})-(\ref{eqn:Ricci curvature for orthogonal and tangent vectors}).
\end{proof}

Finally, we can formulate how the scalar curvature changes while the foliation is warped by a function $f$.

\begin{theo}
The scalar curvature $s^f$ of the warped foliation $M_f$ satisfies
\begin{align}
s^f &= s^{\mathcal{F}} + f^2 (s^{\top})^{\bot} + s^{\bot}
\nonumber\\
&- f^2 \| H^{\mathcal{F}}\|^2 + p\frac{1}{f}\langle H^{\mathcal{F}},\nabla f \rangle - \frac{p(p-1)\|\nabla f\|^2}{f^2} 
\nonumber\\
&+ pf\langle H^{\mathcal{F}},\nabla f\rangle + \sum_{j=1}^{p} (f^2 \|T^{\top} U_j \|^2 - f^2(1-f^2)\|A^{\bot}U_j \|^2
\nonumber\\
&+ (1-f^2) \|S^{\bot} U_j \|^2 + (1-f^2) {\rm tr} ^{\bot}\langle (\nabla_{U_j} S)(\cdot,\cdot), U_j\rangle) - p\frac{{\rm tr}^{\bot} h_f}{f}
\nonumber\\
&- \sum_{i=1}^{q}( (1-f^2)\|A_{X_i}^{\top}\|^2 + 3(1-f^2)\|A^{\bot} X_i\|^2) + \frac{2}{f}\langle H^{\mathcal{F}},\nabla f\rangle 
\nonumber\\
&- p\frac{{\rm tr}^{\bot} h_f}{f} + \frac{1-f^2}{f^2}\langle H^{\bot},H^{\bot}\rangle - \sum_{i=1}^{q}( \frac{1-f^2}{f^2}[\langle S_{X_i}^{\top},A_{X_i}^{\top}\rangle 
\nonumber\\
&- {\rm tr}^{\top} \langle(\nabla_{\cdot} S)(X_i,X_i),\cdot \rangle - \langle S_{X_i}^{\top},(\nabla_{X_i}^{\bot})^{\top}\rangle - \langle S_{X_i}^{\bot},S_{X_i}^{\bot}\rangle]).
\nonumber
\end{align}
where $s^{\top})^{\bot} = \sum_{i=1}^{q} \Ric^{\bot}(U_i,U_i)$.
\end{theo}
\begin{proof}
It follows directly from the formula
\[
s^f (x) = \sum_{i=1}^{m} \ricf (E_i),
\]
and Theorem \ref{thm:Ricci curvature}.
\end{proof}

\section{Examples}

We now will study some examples of warped foliations and its curvatures.

Let $(M,g)$ be a compact 2-dimensional foliated manifold carrying a 1-dimensional foliation $\mathcal{F}$. Let suppose that the sectional curvature of the manifold $\kappa = 0$. Let $(f_n)_{n\in\mathbb{N}}$ be a sequence of constant warping functions such that $f_n=\frac{1}{n}$.

\begin{theo}
 $\lim\limits_{n\to\infty} \kappa^{f_n} = - \langle (\nabla_S)(X,X), U\rangle + \|S(X,U)\|^2$.
\end{theo}
\begin{proof}
By Theorem \ref{thm:sectional curvature},
\begin{align}
\kappa ^f (X,U) &= \kappa (X,U) - \frac{1}{f} h_f(X,X) + 2fXf \langle 
T(U,U), X\rangle
\nonumber\\
&- (1-f^2) [\langle (\nabla_U S)(X,X), U\rangle - \|S(X,U)\|^2]
\nonumber\\
&- f^2(1-f^2) \|A(X,U)\|^2, 
\end{align}
where $X$ and $U$ are vectors orthogonal and tangent to $\mathcal{F}$, respectively. Since $\mathcal{F}$ is a foliation of codimension one, the integrability tensor $A$ of the orthogonal distribution vanishes everywhere. Moreover, $Ef = 0$ for any vector field $E$ on $M$. Again,
\[
\kappa ^{f_n} (X,U) = \kappa (X,U) - (1-f^2) [\langle (\nabla_U S)(X,X), U\rangle - \|S(X,U)\|^2].
\]
Recall that $\kappa(X,U) = \kappa = 0$, and $f_n\to 0$. Finally,
 \[
\lim\limits_{n\to\infty} \kappa^{f_n} = - \langle (\nabla_S)(X,X), U\rangle + \|S(X,U)\|^2.
\]
This ends our proof.
\end{proof}

\begin{coro}
The sectional curvature of a warped by constant functions 1-dimensional Riemmannian foliation on a compact 2-dimensional Riemannian manifold of curvature equal to zero is constant, and remains zero.
\end{coro}
\begin{proof}
 By Theorem \ref{thm:Riemannain foliation and the second fundamental form}, the second fundamental form $S$ vanishes everywhere.
\end{proof}

\end{document}